\documentclass[english]{amsart}
\usepackage[T1]{fontenc}
\usepackage[latin1]{inputenc}
\usepackage{amssymb}

\makeatletter
 \theoremstyle{plain}
\newtheorem{thm}{Theorem}[section]
  \theoremstyle{definition}
  \newtheorem{defn}[thm]{Definition}
  \theoremstyle{plain}
  \newtheorem*{conjecture*}{Conjecture}

\usepackage{babel}
\makeatother
\begin{document}

\title{On The Chaotic Asymptotics of \emph{Ramanujan}'s Entire Function
$A_{q}(z)$ }

\author{Ruiming Zhang}

\date{November 6, 2006}

\email{ruimingzhang@yahoo.com}



\begin{abstract}
We will use a discrete analogue of the classical \emph{Laplace} method
to show that for infinitely many positive integers $n$, the main
term of the asymptotic expansions of scaled confluent basic hypergeometric
functions, including the Ramanujan's entire function $A_{q}(z)$,
could be expressed in $\theta$-functions, and the error term depends
on the ergodic property of certain real scaling parameter. 
\end{abstract}
\maketitle

\section{Introduction}

In \cite{Ismail7} we studied some \emph{Plancherel-Rotach} type asymptotics
for three families of $q$-orthogonal polynomials with logarithmic
scalings. The asymptotics reveals a remarkable pattern which is quite
a contrast to the classical \emph{Plancherel-Rotach} type asymptotics\cite{Deift1,Szego}.
For each family of such polynomials, we found that when a scaling
parameter is above certain critical value, the main term of the asymptotics
involves \emph{Ramanujan}'s entire function $A_{q}(z)$, and the error
term is exponentially small; when the scaling parameter is below the
critical value, if it is rational, then the main term of the asymptotics
has a factor of theta function and the error term is exponentially
small, if the scaling parameter is irrational, then the main term
also contains some theta function, but we could only show that the
error term is no worse than $\mathcal{O}(n^{-1}\log n)$ there. 

In this paper, we will use the same technique to study a \emph{Plancherel-Rotach}
type asymptotics for the confluent basic hypergeometric series, and
we will demonstrate that the asymptotics follows a similar pattern
to the case when a scaling parameter is less than some critical value\cite{Ismail7},
with the only exception that now the base of the theta function explicitly
depends on the type of the confluent basic hypergeometric series.
Furthermore, the error term, in the irrational scaling parameter case,
is closely related to the ergodic property of the number. Generally
speaking, we could improve the error term to $\mathcal{O}(n^{-r})$
, where $r$ is any the positive number less than the chaotic index(we
will define it later in section \ref{sec:Preliminaries}) of the number.
In particular, when the scaling parameter is a real algebraic number
of degree $l$, the error term could not be better than $\mathcal{O}(n^{-r})$
, where the number $r$ is bounded above by $l-1$, and when the scaling
parameter is a \emph{Liouville} number, the error term is better than
$\mathcal{O}(n^{-r})$ for any positive number $r$. 

In section \ref{sec:Preliminaries}, we first introduce the \emph{Ramanujan}'s
entire function $A_{q}(z)$, then define the chaotic index of an irrational
number which is related to its ergodic property. In section \ref{sec:Main-Results},
we will prove asymptotics formulas for $A_{q}(z)$ in detail, and
then just state the results for the confluent basic hypergeometric
series without proofs because the proofs are similar to the corresponding
cases for $A_{q}(z)$.

\section{Preliminaries\label{sec:Preliminaries}}

\subsection{\emph{Ramanujan}'s Entire Function}

Given an arbitrary positive real number $0<q<1$ and an arbitrary
complex number $a$, we define\cite{Gasper,Ismail2} \begin{equation}
(a;q)_{\infty}=\prod_{k=0}^{\infty}(1-aq^{k})\label{eq:2.1}\end{equation}
 and the $q$-shifted factorials \begin{equation}
(a;q)_{n}=\frac{(a;q)_{\infty}}{(aq^{n};q)_{\infty}},\label{eq:2.2}\end{equation}
 for any integer $n\in\mathbb{Z}$. Assuming that $|z|<1$, the $q$-binomial
theorem is\cite{Gasper,Ismail2} \begin{equation}
\frac{(az;q)_{\infty}}{(z;q)_{\infty}}=\sum_{k=0}^{\infty}\frac{(a;q)_{k}}{(q;q)_{k}}z^{k}.\label{eq:2.3}\end{equation}
 Its limiting case \begin{equation}
(z;q)_{\infty}=\sum_{k=0}^{\infty}\frac{q^{k(k-1)/2}}{(q;q)_{k}}(-z)^{k},\label{eq:2.4}\end{equation}
is one of many $q$-exponentials. For any $n\in\mathbb{N}$, let\begin{equation}
R_{1}(a;n):=(aq^{n};q)_{\infty}-1\label{eq:2.5}\end{equation}
 and \begin{equation}
R_{2}(a;n):=\frac{1}{(aq^{n};q)_{\infty}}-1,\label{eq:2.6}\end{equation}
then it is not hard to see that\cite{Ismail7}\begin{equation}
|R_{1}(a;n)|\le\frac{(-aq^{2};q)_{\infty}aq^{n}}{1-q}\label{eq:2.7}\end{equation}
 for $a>0$ and\begin{equation}
|R_{2}(a;n)|\le\frac{aq^{n}}{(1-q)(aq;q)_{\infty}}\label{eq:2.8}\end{equation}
 for $0<aq<1$. 

The \emph{Ramanujan}'s entire function $A_{q}(z)$ is defined as\cite{Ismail2}\begin{equation}
A_{q}(z)=\sum_{k=0}^{\infty}\frac{q^{k^{2}}}{(q;q)_{k}}(-z)^{k}.\label{eq:2.9}\end{equation}
Clearly, it satisfies the following three term recurrence

\begin{equation}
A_{q}(z)-A_{q}(qz)=qzA_{q}(q^{2}z).\label{eq:2.10}\end{equation}
 The function $A_{q}(z)$ appears repeatedly in Ramanujan's work starting
from the Rogers-Ramanujan identities where $A_{q}(-1)$ and $A_{q}(-q)$
are expressed as infinite products\cite{Andrews1}. It is known that
$A_{q}(z)$ has infinitely many zeros, and all of them are positive\cite{Ismail2}.
There have been many researches around $A_{q}(z)$ recently, for example,
to properties of and conjectures about its zeros, \cite{Andrews3},
\cite{Andrews4}, \cite{Hayman}, \cite{Ismail5}. The \emph{Ramanujan}'s
$A_{q}(z)$ is also called as $q$-\emph{Airy} function in the literature,
but this is an unfortunate name, since $A_{q}(z)$ has nothing to
do with the classical \emph{Airy} functions. 

It is clear that \begin{equation}
\left|\frac{(1-q)^{k}}{(q;q)_{k}}q^{k^{2}}(-z)^{k}\right|\le\frac{(q\left|z\right|)^{k}}{k!}\label{eq:2.11}\end{equation}
 for $k=0,1,\dotsc$ and for any complex number $z$. Applying the
\emph{Lebesgue}'s dominated convergent theorem we have\begin{equation}
\lim_{q\to1}A_{q}((1-q)z)=e^{-z}\label{eq:2.12}\end{equation}
 for any $z\in\mathbb{C}$. From (\ref{eq:2.11}) we also have obtained
the inequality \begin{equation}
\left|A_{q}((1-q)z)\right|\le e^{q|z|}\label{eq:2.13}\end{equation}
 for any complex number $z$. 

For any nonzero complex number $z$, we define the theta function
as \begin{equation}
\theta(z;q)=\sum_{k=-\infty}^{\infty}q^{k^{2}/2}z^{k},\label{eq:2.14}\end{equation}
 the \emph{Jacobi}'s triple product formula says that\cite{Gasper,Ismail2}
\begin{equation}
\sum_{k=-\infty}^{\infty}q^{k^{2}/2}z^{k}=(q,-q^{1/2}z,-q^{1/2}/z;q)_{\infty}.\label{eq:2.15}\end{equation}

\subsection{Chaotic Index}

For a positive irrational number $\theta$, \emph{Chebyshev}'s theorem
or its more accurate counterpart \emph{Khinchin}'s theorem implies
that for any real number $\beta$ , there exist infinitely many pairs
of integers $n$ and $m$ with $n>0$ such that\cite{Hua} \begin{equation}
n\theta=m+\beta+\gamma_{n}\quad{\rm with}\quad|\gamma_{n}|\le3/n.\label{eq:2.12}\end{equation}
 Clearly, \emph{Chebyshev}'s theorem and \emph{Khinchin}'s theorem
are saying that the arithmetic progression $\left\{ n\theta\right\} _{n\in\mathbb{Z}}$
is ergodic in $\mathbb{R}$.

\begin{defn}
\label{def:chaotic index}Given real number $\theta$, the \emph{chaotic
index} $\omega(\theta)$ of $\theta$ is defined as the least upper
bound of the set of real numbers $r$ such that for any real number
$\beta$, there exist infinitely many integers $m$ and $n$ with
$n>0$ such that\begin{equation}
|n\theta-\beta-m|\le\frac{C(\theta,r)}{n^{r}},\label{eq:2.13}\end{equation}
 where $C(\theta,r)$ is a positive constant which depends only on
$\theta$ and $r$. 
\end{defn}
\begin{thm}
The chaotic index $\omega(\theta)\ge1$ for any irrational number
$\theta$ . 
\end{thm}
\begin{proof}
This is a direct consequence of the \emph{Chebyshev}'s theorem or
the \emph{Khintchin}'s theorem. 
\end{proof}
Recall that the \emph{irrationality measure} (or \emph{Liouville-Roth
constant}) $\mu(\theta)$ of a real number $\theta$ is defined as
the least upper bound of the set of real numbers $r$ such that\cite{Wikipdedia}
\begin{equation}
0<|n\theta-m|<\frac{1}{n^{r-1}}\label{eq:2.14}\end{equation}
 is satisfied by infinite number of integer pairs $(n,m)$ with $n>0$.
It is clear that we have $\omega(\theta)\le\mu(\theta)-1$, since
given a real number $r<\omega(\theta)$ in the definition\ref{def:chaotic index},
taking $\beta=0$ in (\ref{eq:2.13}), we have \begin{equation}
|n\theta-m|\le\frac{C(\theta,r)}{n^{r}}\label{eq:2.15}\end{equation}
 Recall that a real algebraic number $\theta$ of degree $l$ if it
is a root of an irreducible polynomial of degree $l$ over the integer
ring. \emph{Liouville}'s theorem in number theory says that for a
real algebraic number $\theta$ of degree $l$, there exists a positive
constant $K(\theta)$ such that for any integer $m$ and $n>0$ we
have \begin{equation}
|n\theta-m|>\frac{K(\theta)}{n^{l-1}}.\label{eq:2.16}\end{equation}
 A \emph{Liouville} number is a real number $\theta$ such that for
any positive integer $l$ there exist infinitely many integers $n$
and $m$ with $n>1$ such that\cite{Wikipdedia}\begin{equation}
0<|n\theta-m|<\frac{1}{n^{l-1}},\label{eq:2.17}\end{equation}
 and the terms in the continued fraction expansion of every \emph{Liouville}
number are unbounded. Even though the set of all \emph{Liouville}
numbers is of Lebesgue measure zero, it is known that almost all real
numbers are \emph{Liouville} numbers topologically.

\begin{thm}
The chaotic index $\omega$ has the following properties: 
\begin{enumerate}
\item For any rational number $\theta\in\mathbb{Q}$, its chaotic index
$\omega(\theta)$ is zero; 
\item For any real algebraic number $\theta$ of degree $l$, its chaotic
index $\omega(\theta)\le l-1$; 
\item For any Liouville number $\theta$, its chaotic index $\omega(\theta)=\infty$.\[
\]

\end{enumerate}
\end{thm}
\begin{proof}
Let $\theta=\frac{a}{b}$ for some $a,b\in\mathbb{N}$ with $(a,b)=1$
and $b>0$. If the first claim is not true, then there must be a positive
number $\sigma$, and infinitely many pair of integers $m$ and $n$
with $n>0$ such that \[
|n\theta-m-\beta|\le\frac{C(\theta)}{n^{\sigma}}.\]
 for all $\beta\in\mathbb{R}$. Since $n\theta-m$ are rational numbers
with remainders in the set $\left\{ 0,\frac{1}{b},\dots,\frac{b-1}{b}\right\} $,
so we could pick a real number $\beta$ and find some positive constant
$M(b,\beta)$ such that \[
|n\theta-m-\beta|\ge M(b,\beta)>0,\]
 for any integer pairs $n$ and $m$ with $n>0$, which implies that\[
M(b,\beta)\le\frac{C(\theta)}{n^{\sigma}}.\]
The last inequality is obviously false for $n$ is large enough, which
proves the first assertion. The second and the third assertions follow
from the definition \ref{def:chaotic index}. 
\end{proof}
It is clear that the quadratic irrationals such as $\sqrt{2}$ have
\emph{chaotic index} $\omega=1$. A thorough study on the \emph{chaotic
index} $\omega(\theta)$ as a function in variable $\theta$ is very
desirable, but it is out of the scope of the current paper. It is
also very natural for us to conjecture the following:

\begin{conjecture*}
The chaotic index $\omega(\theta)=\mu(\theta)-1$ for any real number
$\theta$ and $\omega(\theta)=l-1$ if $\theta$ is a real algebraic
number of degree $l$. 
\end{conjecture*}

\section{Main Results\label{sec:Main-Results}}

For any positive real number $t$, we consider the following set \begin{equation}
\mathbb{S}(t)=\left\{ \left\{ nt\right\} :n\in\mathbb{N}\right\} .\label{eq:3.1}\end{equation}
It is clear that $\mathbb{S}(t)\subset[0,1)$ and it is a finite set
when $t$ is rational. In this case, for any $\lambda\in\mathbb{S}(t)$,
there are infinitely many pair of integers $n$ and $m$ with $n>0$
such that \begin{equation}
nt=m+\lambda,\label{eq:3.2}\end{equation}
where\begin{equation}
m=\left\lfloor nt\right\rfloor .\label{eq:3.3}\end{equation}
If $t$ is a positive irrational number, then $\mathbb{S}(t)$ is
a subset of $(0,1)$ with infinite elements, and it is well-known
that $\mathbb{S}(t)$ is uniformly distributed in $(0,1)$. Let $r$
be any positive number less the chaotic index $\omega(t)$, for any
$\beta\in[0,1)$, there are infinitely many pair of integers $n$
and $m$ with $n>0$ such that\begin{equation}
nt=m+\beta+\gamma_{n}\label{eq:3.4}\end{equation}
with\begin{equation}
|\gamma_{n}|\le\frac{C(t,r)}{n^{r}},\label{eq:3.5}\end{equation}
for some constant $C(t,r)$. It is clear that for $n$ large enough,
we may also have\begin{equation}
m=\left\lfloor nt\right\rfloor .\label{eq:3.6}\end{equation}

In the statements of our results, we also need the arithmetic function
$\chi(n)$, the principal character modulo $2$, which is defined
as \begin{equation}
\chi(n)=\begin{cases}
1 & 2\nmid n\\
0 & 2\mid n\end{cases}.\label{eq:3.7}\end{equation}
In the proof of the irrational scaling parameter case, we also make
use of the elementary inequality\begin{equation}
|e^{x}-1|\le|x|e^{|x|}\label{eq:3.8}\end{equation}
 for all $x\in\mathbb{C}$.

\subsection{Asymptotics of $A_{q}(z)$}

For the \emph{Ramanujan}'s entire function $A_{q}(z)$, we have the
following:

\begin{thm}
\label{thm:q_airy}Given an arbitrary positive real number $0<q<1$
and an arbitrary nonzero complex number $u$, we have 
\begin{enumerate}
\item Assuming that $t$ is a positive rational number, for any $\lambda\in\mathbb{S}(t)$,
there are infinitely many pairs of integers $n$ and $m$ satisfying
(\ref{eq:3.2}) and (\ref{eq:3.3}) with $n>0$. For these $n$ and
$m$ we have\begin{equation}
A_{q}(q^{-nt}u)=\frac{(-u)^{\left\lfloor m/2\right\rfloor }\left\{ \theta(-u^{-1}q^{\chi(m)+\lambda};q^{2})+r(n)\right\} }{(q;q)_{\infty}q^{\left\lfloor m/2\right\rfloor (nt-\left\lfloor m/2\right\rfloor )}}\label{eq:3.9}\end{equation}
 with\begin{eqnarray}
|r(n)| & \le & \frac{3(-q^{3};q)_{\infty}\theta(\left|u\right|^{-1}q^{\chi(m)+\lambda};q^{2})}{(q;q)_{\infty}}\nonumber \\
 & \times & \left\{ q^{nt/4}+\frac{q^{m^{2}/16}}{|u|^{\left\lfloor m/4\right\rfloor +1}}\right\} .\label{eq:3.10}\end{eqnarray}

\item Assuming that $t$ is a positive irrational number and $r$ is a positive
number strictly less the chaotic index $\omega(t)$ of $t$, for any
real number $\beta\in[0,1)$, there are infinitely many pair of integers
$n$ and $m$ with $n>0$ satisfying (\ref{eq:3.4}), (\ref{eq:3.5})
and (\ref{eq:3.6}). For such $n$ and $m$ we have\begin{equation}
A_{q}(q^{-nt}u)=\frac{(-u)^{\left\lfloor m/2\right\rfloor }\left\{ \theta(-u^{-1}q^{\chi(m)+\beta};q^{2})+\mathcal{O}\left(n^{-r}\right)\right\} }{(q;q)_{\infty}q^{\left\lfloor m/2\right\rfloor (nt-\left\lfloor m/2\right\rfloor )}}.\label{eq:3.11}\end{equation}
when $n$ is large enough. 
\end{enumerate}
\end{thm}
\begin{proof}
From (\ref{eq:2.9}), we have \begin{equation}
A_{q}(q^{-nt}u)=\sum_{k=0}^{\infty}\frac{q^{k^{2}-knt}}{(q;q)_{k}}(-u)^{k}.\label{eq:3.12}\end{equation}
 In the case that $t$ is a positive rational number, for any $\lambda\in\mathbb{S}(t)$
and $n$, $m$ are large, we have\begin{eqnarray}
A_{q}(q^{-nt}u)(q;q)_{\infty} & = & \sum_{k=0}^{\infty}(q^{k+1};q)_{\infty}q^{k^{2}-km-k\lambda}(-u)^{k}\label{eq:3.13}\\
 & = & s_{1}+s_{2}\nonumber \end{eqnarray}
 where\begin{equation}
s_{1}=\sum_{k=0}^{\left\lfloor m/2\right\rfloor }(q^{k+1};q)_{\infty}q^{k^{2}-km-k\lambda}(-u)^{k}\label{eq:3.14}\end{equation}
 and \begin{equation}
s_{2}=\sum_{k=\left\lfloor m/2\right\rfloor +1}^{\infty}(q^{k+1};q)_{\infty}q^{k^{2}-km-k\lambda}(-u)^{k}.\label{eq:3.15}\end{equation}
 In $s_{1}$ we change the summation from $k$ to $\left\lfloor m/2\right\rfloor -k$,
we have\begin{eqnarray*}
\frac{s_{1}q^{\left\lfloor m/2\right\rfloor (nt-\left\lfloor m/2\right\rfloor )}}{(-u)^{\left\lfloor m/2\right\rfloor }} & = & \sum_{k=0}^{\left\lfloor m/2\right\rfloor }(q^{\left\lfloor m/2\right\rfloor -k+1};q)_{\infty}q^{k^{2}}(-u^{-1}q^{\chi(m)+\lambda})^{k}\\
 & = & \sum_{k=0}^{\infty}q^{k^{2}}(-u^{-1}q^{\chi(m)+\lambda})^{k}\\
 & - & \sum_{k=\left\lfloor m/4\right\rfloor +1}^{\infty}q^{k^{2}}(-u^{-1}q^{\chi(m)+\lambda})^{k}\\
 & + & \sum_{k=0}^{\left\lfloor m/4\right\rfloor }q^{k^{2}}(-u^{-1}q^{\chi(m)+\lambda})^{k}\left((q^{\left\lfloor m/2\right\rfloor -k+1};q)_{\infty}-1\right)\\
 & + & \sum_{k=\left\lfloor m/4\right\rfloor +1}^{\left\lfloor m/2\right\rfloor }(q^{\left\lfloor m/2\right\rfloor -k+1};q)_{\infty}q^{k^{2}}(-u^{-1}q^{\chi(m)+\lambda})^{k}\\
 & = & \sum_{k=0}^{\infty}q^{k^{2}}(-u^{-1}q^{\chi(m)+\lambda})^{k}+s_{11}+s_{12}+s_{13}.\end{eqnarray*}
 Clearly,\begin{eqnarray*}
|s_{11}+s_{13}| & \le & 2\sum_{k=\left\lfloor m/4\right\rfloor +1}^{\infty}q^{k^{2}}(\left|u\right|^{-1}q^{\chi(m)+\lambda})^{k}\\
 & \le & \frac{2q^{m^{2}/16}}{|u|^{\left\lfloor m/4\right\rfloor +1}}\theta(\left|u\right|^{-1}q^{\chi(m)+\lambda};q^{2}).\end{eqnarray*}
 By (\ref{eq:2.7}), for $0\le k\le\left\lfloor m/4\right\rfloor $,
we have\begin{equation}
\left|(q^{\left\lfloor m/2\right\rfloor -k+1};q)_{\infty}-1\right|=\frac{(-q^{3};q)_{\infty}}{(q;q)_{\infty}}q^{nt/4},\label{eq:3.16}\end{equation}
 then\begin{eqnarray*}
|s_{12}| & \le & \frac{(-q^{3};q)_{\infty}q^{nt/4}}{(q;q)_{\infty}}\sum_{k=0}^{\infty}q^{k^{2}}(\left|u\right|^{-1}q^{\chi(m)+\lambda})^{k}\\
 & \le & \frac{q^{nt/4}(-q^{3};q)_{\infty}\theta(\left|u\right|^{-1}q^{\chi(m)+\lambda};q^{2})}{(q;q)_{\infty}},\end{eqnarray*}
 hence\begin{equation}
\frac{s_{1}q^{\left\lfloor m/2\right\rfloor (nt-\left\lfloor m/2\right\rfloor )}}{(-u)^{\left\lfloor m/2\right\rfloor }}=\sum_{k=0}^{\infty}q^{k^{2}}(-u^{-1}q^{\chi(m)+\lambda})^{k}+r_{1}(n)\label{eq:3.17}\end{equation}
 with\begin{eqnarray}
|r_{1}(n)| & \le & \frac{2(-q^{3};q)_{\infty}\theta(\left|u\right|^{-1}q^{\chi(m)+\lambda};q^{2})}{1-q}\label{eq:3.18}\\
 & \times & \left\{ q^{nt/4}+\frac{q^{m^{2}/16}}{|u|^{\left\lfloor m/4\right\rfloor +1}}\right\} .\nonumber \end{eqnarray}
 In $s_{2}$ we change the summation from $k$ to $k+\left\lfloor m/2\right\rfloor $\begin{eqnarray*}
\frac{s_{2}q^{\left\lfloor m/2\right\rfloor (nt-\left\lfloor m/2\right\rfloor )}}{(-u)^{\left\lfloor m/2\right\rfloor }} & = & \sum_{k=1}^{\infty}(q^{\left\lfloor m/2\right\rfloor +k+1};q)_{\infty}q^{k^{2}}(-uq^{-\chi(m)-\lambda})^{k}\\
 & = & \sum_{k=1}^{\infty}q^{k^{2}}(-uq^{-\chi(m)-\lambda})^{k}\\
 & + & \sum_{k=1}^{\infty}q^{k^{2}}(-uq^{-\chi(m)-\lambda})^{k}\left[(q^{\left\lfloor m/2\right\rfloor +1};q)_{\infty}-1\right]\\
 & = & \sum_{k=-\infty}^{-1}q^{k^{2}}(-u^{-1}q^{\chi(m)+\lambda})^{k}+r_{2}(n).\end{eqnarray*}
 By (\ref{eq:2.7}), for $k\ge1$\begin{equation}
\left|(q^{\left\lfloor m/2\right\rfloor +1};q)_{\infty}-1\right|\le\frac{q^{nt/2}(-q^{3};q)_{\infty}}{(q;q)_{\infty}},\label{eq:3.19}\end{equation}
 then\begin{equation}
|r_{2}(n)|\le\frac{q^{nt/2}(-q^{3};q)_{\infty}\theta(\left|u\right|^{-1}q^{\chi(m)+\lambda};q^{2})}{(q;q)_{\infty}}.\label{eq:3.20}\end{equation}
 Thus we have proved that\begin{equation}
A_{q}(q^{-nt}u)=\frac{(-u)^{\left\lfloor m/2\right\rfloor }\left\{ \theta(-u^{-1}q^{\chi(m)+\lambda};q^{2})+r(n)\right\} }{(q;q)_{\infty}q^{\left\lfloor m/2\right\rfloor (nt-\left\lfloor m/2\right\rfloor )}}\label{eq:3.21}\end{equation}
 with\begin{eqnarray}
|r(n)| & \le & \frac{3(-q^{3};q)_{\infty}\theta(\left|u\right|^{-1}q^{\chi(m)+\lambda};q)}{1-q}\label{eq:3.22}\\
 & \times & \left\{ q^{nt/4}+\frac{q^{m^{2}/16}}{|u|^{\left\lfloor m/4\right\rfloor +1}}\right\} .\nonumber \end{eqnarray}
 for $n$, $m$ and $\lambda$ satisfying (\ref{eq:3.2}) , (\ref{eq:3.3})with
$n$ and $m$ are large enough.

In the case that $t$ is a positive irrational number, for any real
number $\beta\in[0,1)$, when $n$ and $m$ satisfy (\ref{eq:3.4}),
(\ref{eq:3.5}) and (\ref{eq:3.6}) and $n$ is large enough, we take\begin{equation}
\nu_{n}=\left\lfloor -\frac{\log n}{\log q}\right\rfloor ,\label{eq:2.23}\end{equation}
 then we have \begin{eqnarray}
A_{q}(q^{-nt}u)(q;q)_{\infty} & = & \sum_{k=0}^{\infty}(q^{k+1};q)_{\infty}q^{k^{2}-km-k\beta-k\gamma_{n}}(-u)^{k}\label{eq:3.24}\\
 & = & s_{1}+s_{2},\nonumber \end{eqnarray}
 with\begin{equation}
s_{1}=\sum_{k=0}^{\left\lfloor m/2\right\rfloor }(q^{k+1};q)_{\infty}q^{k^{2}-km-k\beta-k\gamma_{n}}(-u)^{k}\label{eq:3.25}\end{equation}
 and\begin{equation}
s_{2}=\sum_{k=\left\lfloor m/2\right\rfloor +1}^{\infty}(q^{k+1};q)_{\infty}q^{k^{2}-km-k\beta-k\gamma_{n}}(-u)^{k}.\label{eq:3.26}\end{equation}
 In $s_{1}$, we change summation from $k$ to $\left\lfloor m/2\right\rfloor -k$,
\begin{eqnarray*}
\frac{s_{1}q^{\left\lfloor m/2\right\rfloor (nt-\left\lfloor m/2\right\rfloor )}}{(-u)^{\left\lfloor m/2\right\rfloor }} & = & \sum_{k=0}^{\left\lfloor m/2\right\rfloor }(q^{\left\lfloor m/2\right\rfloor -k+1};q)_{\infty}q^{k^{2}}(-u^{-1}q^{\chi(m)+\beta+\gamma_{n}})^{k}\\
 & = & \sum_{k=0}^{\infty}q^{k^{2}}(-u^{-1}q^{\chi(m)+\beta})^{k}\\
 & - & \sum_{k=\nu_{n}+1}^{\infty}q^{k^{2}}(-u^{-1}q^{\chi(m)+\beta})^{k}\\
 & + & \sum_{k=0}^{\nu_{n}}q^{k^{2}}(-u^{-1}q^{\chi(m)+\beta})^{k}\left(q^{k\gamma_{n}}-1\right)\\
 & + & \sum_{k=0}^{\nu_{n}}q^{k^{2}}(-u^{-1}q^{\chi(m)+\beta+\gamma_{n}})^{k}\left\{ (q^{\left\lfloor m/2\right\rfloor -k+1};q)_{\infty}-1\right\} \\
 & + & \sum_{k=\nu_{n}+1}^{\left\lfloor m/2\right\rfloor }q^{k^{2}}(-u^{-1}q^{\chi(m)+\beta+\gamma_{n}})^{k}(q^{\left\lfloor m/2\right\rfloor -k+1};q)_{\infty}\\
 & = & \sum_{k=0}^{\infty}q^{k^{2}}(-u^{-1}q^{\chi(m)+\beta})^{k}+s_{11}+s_{12}+s_{13}+s_{14},\end{eqnarray*}
 thus there exists some constant $c_{11}(q,u)$ such that\begin{eqnarray*}
|s_{11}+s_{14}| & \le & 2\sum_{k=\nu_{n}+1}^{\infty}q^{k^{2}}(|u|^{-1}q^{\chi(m)+\beta})^{k}\\
 & \le & c_{11}(q,u)n^{-r}\log n\end{eqnarray*}
 for $n$ large enough, and there exists some constant $c_{12}(q,u)$
such that\begin{eqnarray*}
|s_{12}+s_{13}| & \le & c_{12}(q,u)n^{-r}\log n,\end{eqnarray*}
 for $n$ large enough. Thus for $n$, $m$ and $\beta$ satisfying
(\ref{eq:3.4}) , (\ref{eq:3.5}) and (\ref{eq:3.6}) and $n$ is
large enough, there exists some constant $c_{1}(q,u)$ such that\begin{equation}
\frac{s_{1}q^{\left\lfloor m/2\right\rfloor (nt-\left\lfloor m/2\right\rfloor )}}{(-u)^{\left\lfloor m/2\right\rfloor }}=\sum_{k=0}^{\infty}q^{k^{2}}(-u^{-1}q^{\chi(m)+\beta})^{k}+e_{1}(n)\label{eq:3.27}\end{equation}
 with\begin{equation}
|e_{1}(n)|\le c_{1}(q,u)n^{-r}\log n.\label{eq:3.28}\end{equation}
 Similarly, in $s_{2}$ we change summation from $k$ to $k+\left\lfloor m/2\right\rfloor $,
and we could show that there for $n$, $m$ and $\beta$ satisfying
(\ref{eq:3.4}) , (\ref{eq:3.5}) and (\ref{eq:3.6}) and $n$ is
large enough, there exists some constant $c_{2}(q,u)$ such that\begin{equation}
\frac{s_{2}q^{\left\lfloor m/2\right\rfloor (nt-\left\lfloor m/2\right\rfloor )}}{(-u)^{\left\lfloor m/2\right\rfloor }}=\sum_{k=-\infty}^{-1}q^{k^{2}}(-u^{-1}q^{\chi(m)+\beta})^{k}+e_{2}(n)\label{eq:3.29}\end{equation}
 with\begin{equation}
|e_{2}(n)|\le c_{2}(q,u)n^{-r}\log n.\label{eq:3.30}\end{equation}
 Thus we have\begin{equation}
A_{q}(q^{-nt}u)=\frac{(-u)^{\left\lfloor m/2\right\rfloor }\left\{ \theta(-u^{-1}q^{\chi(m)+\beta};q^{2})+\mathcal{O}\left(n^{-r}\log n\right)\right\} }{(q;q)_{\infty}q^{\left\lfloor m/2\right\rfloor (nt-\left\lfloor m/2\right\rfloor )}}\label{eq:3.31}\end{equation}
 for $n$, $m$ and $\beta$ satisfying (\ref{eq:3.4}) , (\ref{eq:3.5})
and (\ref{eq:3.6}) and $n$ is large enough. Since $\log n=\mathcal{O}(n^{-\epsilon})$
for any positive number $\epsilon$, and any $r$ is an arbitrary
positive number strictly less the chaotic index $\omega(t)$, thus
we have proved the last assertion.
\end{proof}

\subsection{Asymptotics of Confluent Basic Hypergeometric Functions}

The phenomenon demonstrated with $A_{q}(z)$ in theorem\ref{thm:q_airy}
is universal for a class of entire basic hypergeometric function of
type\begin{equation}
f(z)=\sum_{k=0}^{\infty}\frac{(a_{1},\dotsc,a_{r};q)_{k}q^{lk^{2}}}{(b_{1},\dotsc,b_{s};q)_{k}}z^{k},\label{eq:3.32}\end{equation}
 with $l>0$, where \begin{equation}
(a_{1},\dotsc,a_{r};q)_{k}=\prod_{j=1}^{r}(a_{j};q)_{k}.\label{eq:3.33}\end{equation}
 A basic hypergeometric series is formally defined as\begin{equation}
_{r}\phi_{s}\left(\begin{array}{c|c}
\begin{array}{c}
a_{1},\dotsc,a_{r}\\
b_{1},\dots,b_{s}\end{array} & q,z\end{array}\right)=\sum_{k=0}^{\infty}\frac{(a_{1},\dotsc,a_{r};q)_{k}}{(b_{1},\dotsc,b_{s};q)_{k}}z^{k}\left(-q^{(k-1)/2}\right)^{k(s+1-r)},\label{eq:3.34}\end{equation}
 when $s+1-r>0$, it is clear that the function\begin{equation}
_{r}\phi_{s}\left(\begin{array}{c|c}
\begin{array}{c}
a_{1},\dotsc,a_{r}\\
b_{1},\dots,b_{s}\end{array} & q,-zq^{(s+1-r)/2}\end{array}\right)\label{eq:3.35}\end{equation}
 is of the form (\ref{eq:3.32}) with \begin{equation}
l=\frac{s+1-r}{2}.\label{eq:3.36}\end{equation}
 For a clean statement of the following theorem, we also define\begin{equation}
c(r,s;q):=\frac{(a_{1},\dotsc,a_{r};q)_{\infty}}{(b_{1},\dotsc,b_{s};q)_{\infty}}.\label{eq:3.37}\end{equation}
 We have similar results for (\ref{eq:3.32}) as theorem \ref{thm:q_airy}
for $A_{q}(z)$, and their proof is very similar to the proof of theorem
\ref{thm:q_airy} without any significant changes.

\begin{thm}
\label{thm:confluent}Assume that \begin{equation}
0\le a_{1},\dotsc,a_{r},b_{1},\dotsc,bs<1,\label{eq:3.38}\end{equation}
 and $0<q<1$. For any nonzero complex number $u$, we have 
\begin{enumerate}
\item Given any positive rational number $t$, let $\lambda\in\mathbb{S}(t)$,
there are infinitely many pair of integers $n$ and $m$ satisfying
(\ref{eq:3.2}) and (\ref{eq:3.3}) with $n>0$. For these $n$ and
$m$ we have\begin{equation}
f(q^{-lnt}u)=\frac{(-u)^{\left\lfloor m/2\right\rfloor }\left\{ \theta(u^{-1}q^{l\chi(m)+l\lambda};q^{2l})+r(n)\right\} }{c(s,r;q)q^{l\left\lfloor m/2\right\rfloor (nt-\left\lfloor m/2\right\rfloor )}}\label{eq:3.39}\end{equation}
 with\begin{eqnarray}
|r(n)| & \le & \frac{2^{r+s+2}\theta(\left|u\right|^{-1}q^{l\chi(m)+l\lambda};q^{2l})}{(q;q)_{\infty}^{r+s}}\nonumber \\
 & \times & \left\{ \frac{{\displaystyle \prod_{j=1}^{s}(-b_{j}q^{2};q)_{\infty}}}{{\displaystyle \prod_{j=1}^{r}}(a_{j}q;q)_{\infty}}q^{nt/4}+\frac{q^{lm^{2}/16}}{|u|^{\left\lfloor m/4\right\rfloor +1}}\right\} .\label{eq:3.40}\end{eqnarray}
 for $n$, $m$ and $\lambda$ satisfying (\ref{eq:3.2}) and (\ref{eq:3.3})
with $n$ is large enough. 
\item Given a positive irrational number $t$ and a positive number $r<\omega(t)$,
for any real number $\beta\in[0,1)$, there are infinitely many pair
of integers $n$ and $m$ satisfying (\ref{eq:3.4}), (\ref{eq:3.5})
and (\ref{eq:3.6}) with $n>0$. For such $n$ and $m$ we have\begin{equation}
f(q^{-lnt}u)=\frac{(-u)^{\left\lfloor m/2\right\rfloor }\left\{ \theta(u^{-1}q^{l\chi(m)+l\lambda};q^{2l})+\mathcal{O}\left(n^{-r}\right)\right\} }{c(s,r;q)q^{l\left\lfloor m/2\right\rfloor (nt-\left\lfloor m/2\right\rfloor )}}.\label{eq:3.41}\end{equation}
 for $n$ is large enough. 
\end{enumerate}
\end{thm}

\end{document}